\RequirePackage{ifpdf}
\ifpdf 
\documentclass[pdftex]{sigma}
\else
\documentclass{sigma}
\fi

\usepackage{mathrsfs}

\newcommand{\p}{\mathbb{P}}
\newcommand{\e}{\mathbb{E}}
\newcommand{\re}{\mathbb{R}}

\newcommand{\fc}{{}_1\mathscr{F}_1}

\begin{document}
\allowdisplaybreaks

\renewcommand{\thefootnote}{$\star$}

\renewcommand{\PaperNumber}{074}

\FirstPageHeading

\ShortArticleName{First Hitting Time of the Boundary of the Weyl Chamber by Radial Dunkl Processes}

\ArticleName{First Hitting Time of the Boundary \\
of the Weyl Chamber by Radial Dunkl Processes\footnote{This paper is a contribution to the Special
Issue on Dunkl Operators and Related Topics. The full collection
is available at
\href{http://www.emis.de/journals/SIGMA/Dunkl_operators.html}{http://www.emis.de/journals/SIGMA/Dunkl\_{}operators.html}}}

\Author{Nizar DEMNI}

\AuthorNameForHeading{N. Demni}

\Address{SFB 701, Fakult\"at f\"ur Mathematik, Universit\"at Bielefeld, Deutschland}
\Email{\href{mailto:demni@math.uni-bielefeld.de}{demni@math.uni-bielefeld.de}}

\ArticleDates{Received July 01, 2008, in f\/inal form October 24,
2008; Published online November 04, 2008}

\Abstract{We provide two equivalent approaches for computing the tail distribution of the f\/irst hitting time of the boundary of the Weyl chamber by a radial Dunkl process. The f\/irst approach is based on a spectral problem with initial value. The second one expresses the tail distribution by means of the $W$-invariant Dunkl--Hermite polynomials.
Illustrative examples are given by the irreducible root systems of types $A$, $B$, $D$. The paper ends with an interest in the case of Brownian motions for which our formulae take determinantal forms.}


\Classification{33C20; 33C52; 60J60; 60J65}

\renewcommand{\thefootnote}{\arabic{footnote}}
\setcounter{footnote}{0}

\section{Motivation}
The f\/irst exit time from cones by a multidimensional Brownian motion has been of great interest for mathematicians~\cite{Ban,Dou} and for theoretical physicists as well~\cite{Com}. An old and famous example of cones is provided by root systems in a f\/inite dimensional Euclidean space, say $(V,\langle\cdot\rangle)$~\cite{Hum}. More precisely, a root system $R$ in $V$ is a collection of non zero vectors from $V$ such that $\sigma_{\alpha}(R) = R$ for all $\alpha \in R$, where $\sigma_{\alpha}$ is the ref\/lection with respect to the hyperplane orthogonal to $\alpha$:
\begin{gather*}
\sigma_{\alpha} (x) =  x - 2\frac{\langle \alpha, x\rangle}{\langle\alpha,\alpha\rangle} \alpha , \qquad x \in V.
 \end{gather*}
A simple system $S$ is  a basis of $\textrm{span}(R)$ which induces a total ordering in $R$. A  root $\alpha$ is positive if it is a positive linear combination of elements of~$S$. The set of positive roots is called a positive subsystem and is denoted by $R_+$. The cone $C$ associated with $R$, known as the positive Weyl chamber,  is def\/ined by
\begin{gather*}
C := \{x \in V, \  \langle \alpha , x \rangle > 0 \  \forall \, \alpha \in R_+\} = \{x \in V, \  \langle \alpha , x \rangle > 0 \  \forall \, \alpha \in S\}.
\end{gather*}
For such cones, explicit formulae for the f\/irst exit time were given in \cite{Dou} and involve Pfaf\/f\/ians of skew-symmetric matrices while \cite{Ban} covers more general cones.
During the last decade, a~dif\/fusion process valued in $\overline{C}$, the topological closure of $C$, was introduced and studied in a series of papers~\cite{Chy,Gal1,Gal2,Gal3,GG1} and generalizes the ref\/lected Brownian motion, that is, the absolute value of a real Brownian motion. This dif\/fusion, known as the radial Dunkl process, is associated with a root system and depends on a set of positive parameters called a multiplicity function often denoted $k$. The latter is def\/ined on the orbits of the action of the group generated by all the ref\/lections $\sigma_{\alpha}$, $\alpha \in R$, the ref\/lection group, and is constant on each orbit. Let $X$ denote a~radial Dunkl process starting at $x \in C$ and def\/ine the f\/irst hitting time of $\partial C$ by
\begin{gather*}
T_0 := \inf\{t, X_t \in \partial C\}.
\end{gather*}
Let $l: \alpha \in R \mapsto  k(\alpha)-1/2$ be the so called index function~\cite{Chy}. Then it was shown in \cite{Chy} that $T_0 < \infty$ almost surely (hereafter a.s.)  if $-1/2 \leq l(\alpha) < 0$ for at least one $\alpha \in R$, and $T_0 = \infty$ a.s.\ if $l(\alpha) \geq 0$ for all $\alpha \in R$. Moreover, the tail distribution of $T_0$ may be computed from the absolute-continuity relations derived in \cite{Chy}. Two cases are distinguished
\begin{itemize}\itemsep=0pt
\item $0 \leq l(\alpha) \leq 1/2$ for all $\alpha \in R$ with at least one $\alpha$ such that $0 < l(\alpha) \leq1/2$: the radial Dunkl process with index $-l$ hits $\partial C$ a.s.;
\item $-1/2 \leq l(\alpha) < 0$ for at least one $\alpha$ and $l(\beta) \geq  0$ for at least one $\beta$: $X$ itself hits $\partial C$ a.s.
\end{itemize}
In the f\/irst case, we shall address the problem in its whole generality and specialize our results to the types $A$, $B$, $D$, while in the second case we shall restrict ourselves to the type~$B$.
One of the reasons is that the second case needs two values of the multiplicity function or equivalently two orbits. Another reason is that, after a suitable change of the index function, we are led to the f\/irst case, that is, to the case when both indices are positive.

After this panorama, we present another approach which is equivalent to the one used before and has the merit to express the tail distribution through the $W$-invariant parts of the so-called Dunkl--Hermite polynomials\footnote{They were called generalized Hermite polynomials in \cite{Ros1} but we prefer calling them as above to avoid the confusion with the generalized Hermite polynomials introduced by Lassalle~\cite{Lass2}.}. The last part is concerned with determinantal formulae obtained for the f\/irst hitting time of $\partial C$ by a multi-dimensional Brownian motion. These formulae have to be compared to the ones obtained in \cite{Dou}.

\section{First approach}
\subsection{First case} \label{FF}
Let us denote by $\p_x^{l}$ the law of $(X_t)_{t \geq 0}$ starting from $x \in C$ and of index function $l$. Let $\e_x^{l}$ be the corresponding expectation. Recall that~\cite[p.~38, Proposition 2.15c]{Chy}, if $l(\alpha) \geq 0$ for all $\alpha \in R_+$, then
\begin{gather*}
\p_x^{-l}(T_0 > t) = \e_x^l\left[\left(\prod_{\alpha \in R_+}\frac{\langle \alpha,X_t \rangle}{\langle \alpha,x\rangle }\right)^{-2l(\alpha)}\right].
\end{gather*}
Recall that the semi group density of $X$ is given by \cite{Chy,Ros0}
\begin{gather*}
p_t^k(x,y) =\frac{1}{c_kt^{\gamma + m/2}} e^{-(|x|^2 + |y|^2)/2t}D_k^W\left(\frac{x}{\sqrt t},\frac{y}{\sqrt t}\right) \omega_k^2(y), \qquad x,y \in \overline{C},
\end{gather*}
where $\gamma = \sum\limits_{\alpha \in R_+}k(\alpha)$ and $m= \dim V$ is the rank of $R$. The weight function $\omega_k$ is given by
\begin{gather*}
\omega_k(y) = \prod_{\alpha \in R_+} \langle \alpha,y \rangle ^{k(\alpha)}
\end{gather*}
and $D_k^W$ is the generalized Bessel function. Thus,
\begin{align*}
\p_x^{-l}(T_0 > t) & = \prod_{\alpha \in R_+}\langle \alpha,x\rangle ^{2l(\alpha)}\frac{e^{-|x|^2/2t}}{c_kt^{\gamma + m/2}}\int_{C}
e^{-|y|^2/2t}D_k^W\left(\frac{x}{\sqrt t}, \frac{y}{\sqrt t}\right)\prod_{\alpha \in R_+}\langle \alpha,y\rangle dy
\\& = \prod_{\alpha \in R_+}\langle \alpha,\frac{x}{\sqrt t} \rangle ^{2l(\alpha)}\frac{e^{-|x|^2/2t}}{c_k}\int_{C}e^{-|y|^2/2}
D_k^W\left(\frac{x}{\sqrt t}, y\right)\prod_{\alpha \in R_+}\langle \alpha,y\rangle dy
\\&  =  \prod_{\alpha \in R_+}\langle \alpha,\frac{x}{\sqrt t}\rangle ^{2l(\alpha)}\frac{e^{-|x|^2/2t}}{c_k} g\left(\frac{x}{\sqrt t}\right),
\end{align*}
where
\begin{gather}\label{E1}
g(x) := \int_{C}e^{-|y|^2/2}D_k^W(x,y)\prod_{\alpha \in R_+}\langle \alpha,y\rangle dy.
\end{gather}
Our key result is stated as follows:
\begin{theorem}\label{T1}
Let $T_i$ be the $i$-th Dunkl derivative and $\Delta_k = \sum\limits_{i=1}^m T_i^2$ the Dunkl Laplacian. Define
\begin{gather*}
\mathscr{J}_k^x :=  -\Delta_k^x + \sum_{i =1}^m x_i\partial_i^x : = -\Delta_k^x + E_1^x,
\end{gather*}
where $E_1^x := \sum\limits_{i=1}^mx_i\partial_i^x$ is the Euler operator and the superscript indicates the derivative action. Then
\begin{gather*}
\mathscr{J}_k^x \big[e^{-|y|^2/2}D_k^W(x, y)\big] = E_1^y\big[e^{-|y|^2/2}D_k^W(x, y)\big].
\end{gather*}
\end{theorem}
\begin{proof} Recall that if $f$ is $W$-invariant then $T_i^xf = \partial_i^xf$ and that $T_i^xD_k(x,y) = y_iD_k(x,y)$ (see~\cite{Ros1}). On the one hand~\cite{Ros1}
\begin{align*}
\Delta_k^x D_k^W(x, y) = \sum_{i=1}^my_i^2\sum_{w \in W}D_k(x,wy) = |y|^2D_k^W(x,y).
\end{align*}
On the other hand,	
\begin{align*}
E_1^x D_k^W(x, y) & = \sum_{w \in W}\sum_{i=1}^mx_i T_i^xD_k(x,wy) = \sum_{w \in W}\sum_{i=1}^m(x_i)(wy)_i D_k(x,wy)
\\& = \sum_{w \in W}\langle x,wy\rangle  D_k(x,wy) = \sum_{w \in W}\langle w^{-1}x,y\rangle D_k(x,wy)
\\& = E_1^y D_k^W(x,y),
\end{align*}
where the last equality follows from  $D_k(x,wy) = D_k(w^{-1}x,y)$ since $D_k(wx,wy) = D_k(x,y)$ for all $w \in W$. The result follows from an easy computation.
\end{proof}
\begin{remark}
The appearance of the operator $\mathscr{J}_k$ is not a mere coincidence and will be explained when presenting the second approach.
\end{remark}
 \begin{corollary}
$g$ is an eigenfunction of $-\mathscr{J}_k$ corresponding to the eigenvalue $m+|R_+|$.
\end{corollary}
\begin{proof} Theorem~\ref{T1} and an integration by parts give
\begin{align*}
-\mathscr{J}_k^x g(x) &= - \int_{C}E_1^y\left[e^{-|y|^2/2}D_k^W(x, y)\right]\prod_{\alpha \in R_+}\langle \alpha,y\rangle dy
\\& = -\sum_{i=1}^m \int_{C}y_i\prod_{\alpha \in R_+}\langle \alpha,y \rangle \partial_i^y\left[e^{-|y|^2/2}D_k^W(x, y)\right] \,dy
\\& = \sum_{i=1}^m \int_C e^{-|y|^2/2}D_k^W(x, y) \partial_i\left[y_i\prod_{\alpha \in R_+}\langle \alpha,y\rangle \right] dy
\\& = \int_C  e^{-|y|^2/2}D_k^W(x, y)\prod_{\alpha \in R_+}\langle \alpha,y\rangle \sum_{i=1}^m \left[1 + \sum_{\alpha \in R_+}\frac{\alpha_iy_i}{\langle\alpha,y\rangle}\right] dy.
\end{align*}
The proof ends after summing over $i$.
\end{proof}

\begin{remark}
The crucial advantage in our approach is that we do not need the explicit expression of $D_k^W$. Moreover, though $D_k^W$ can be expressed through multivariate hypergeometric series~\cite{Bak,Dem,Dem1}, it cannot in general help to compute the function~$g$.
\end{remark}

\subsubsection[The $A$-type]{The $\boldsymbol{A}$-type}

This root system is characterized by
\begin{alignat*}{3}
&R  =  \{\pm (e_i - e_j),\  1 \leq i < j \leq m\},\qquad  &&   R_+  =  \{e_i - e_j,\  1 \leq i < j \leq m\},&\\
& S  =  \{e_i - e_{i+1}, \  1 \leq i \leq m-1\},\qquad & & C = \{y \in \re^m, \  y_1 > y_2 > \cdots > y_m\}.&
\end{alignat*}
The ref\/lection group $W $ is the permutations group $S_m$ and there is one orbit so that $k = k_1 > 0$ thereby $\gamma = k_1m(m-1)/2$.
Moreover, $D_k^W$ is given by (see~\cite[p.~212--214]{Bak})\footnote{Authors used another normalization so that there is factor $\sqrt 2$ in both arguments. Moreover the Jack para\-me\-ter denoted there by $\alpha$ is the inverse of $k_1$.}
\begin{gather*}
\frac{1}{|W|}D_k^W(x,y) =  {}_0F_0^{(1/k_1)}(x,y).
\end{gather*}
Hence, letting $y \mapsto V(y)$ be the Vandermonde function, one writes:
\begin{align*}
\p_x^{-l}(T_0 > t)  &=  V\left(\frac{x}{\sqrt{t}}\right)^{2l_1}\frac{|W|e^{-|x|^2/2t}}{c_k} \int_{C}e^{-|y|^2/2} {}_0F_0^{(1/k_1)}\left(\frac{x}{\sqrt t}, y\right) V(y) dy,
\end{align*}
where $l_1 = k_1 -1/2$. Besides, $\mathscr{J}_k$ acts on $W$-invariant functions as
\begin{gather*}
-\mathscr{J}_k^x = D_0^x - E_1^x := \sum_{i=1}^m \partial_i^{2,x} + 2k_1\sum_{i\neq j}\frac{1}{x_i - x_j}\partial_i^x -  \sum_{i=1}^mx_i \partial_i^x.
\end{gather*}
Finally, since $g$ is $W$-invariant, then
\begin{gather*}
\left(D_0^x - E_1^x\right)g(x)   =   m\frac{m+1}{2}g(x), \\
g(0) = m!\int_C e^{-|y|^2/2} V(y) dy   =    \int_{\re^m} e^{-|y|^2/2} |V(y)| dy.
\end{gather*}
In order to write down $g$, let us recall that the multivariate Gauss hypergeometric function ${}_2F_1^{(1/k_1)}(e,b,c,\cdot)$ (see~\cite{Bee} for the def\/inition) is the unique symmetric eigenfunction that equals $1$ at $0$ of (see \cite[p.~585]{Bee})
\begin{gather}
\sum_{i=1}^mz_i(1-z_i)\partial_i^{2,z} + 2k_1\sum_{i\neq j}\frac{z_i(1-z_i)}{z_i - z_j} \partial_i^z\nonumber\\
\qquad{}+ \sum_{i=1}^m\left[c- k_1(m-1) -\left(e+b+1 - k_1(m-1)\right) z_i\right]\partial_i^z\label{GF}
\end{gather}
associated with the eigenvalue $meb$. Letting $z = (1/2)(1- x/\sqrt{b})$, that is $z_i = (1/2)(1- x_i/\sqrt{b})$ for all $1 \leq i \leq m$ and $e=(m+1)/2$,
\begin{gather*}
c = k_1(m-1) + \frac{1}{2}[e+b+1 - k_1(m-1)] = \frac{b}{2} + \frac{k_1}{2}(m-1) + \frac{m+3}{4},
\end{gather*}
then, the resulting function is an eigenfunction of
\begin{gather*}
\sum_{i=1}^m \left(1-\frac{x_i^2}{b}\right)\partial_i^{2,x} + 2k_1\sum_{i\neq j}\frac{(1-x_i^2/b)}{x_i - x_j} \partial_i^x - \sum_{i=1}^m \left(b+\frac{m+3}{2} - k_1(m-1)\right)\frac{x_i}{b}\partial_i^x
\end{gather*}
and $D_0^x - E_1^x$ is the limiting operator as $b$ tends to inf\/inity. Hence,
\begin{proposition}For $1/2 < k_1 \leq 1$,
\begin{gather*}
g(x) = g(0) C(m,k_1)\lim_{b \rightarrow \infty} {}_2F_1^{(1/k_1)}\left[\frac{m+1}{2}, b, \frac{b}{2} + \frac{k_1}{2}(m-1) + \frac{m+3}{4}, \frac{1}{2}\left(1-\frac{x}{\sqrt{b}}\right)\right],
\end{gather*}
where
\begin{gather*}
C(m, k_1)^{-1} = \lim_{b \rightarrow \infty}{}_2F_1^{(1/k_1)}\left(\frac{m+1}{2}, b, \frac{b}{2} + \frac{k_1}{2}(m-1) + \frac{m+3}{4}, \frac{1}{2}\right).
\end{gather*}\end{proposition}

\begin{remark}
One cannot exchange the inf\/inite sum and the limit operation. Indeed, expand the generalized Pochhammer symbol as (see \cite[p.~191]{Bak})
\begin{gather*}
(a)_{\tau} = \prod_{i=1}^m (a - k_1(i-1))_{\tau_i} = \prod_{i=1}^m\frac{\Gamma(a - k_1(i-1) + \tau_i)}{\Gamma(a - k_1(i-1))}
\end{gather*}
and use Stirling formula to see that each term in the above product is equivalent to
\[
(a + k_1(m-1)+\tau_i)^{\tau_i}
 \]
 for large enough positive $a$. Moreover, since $J_{\tau}^{(1/k_1)}$ is homogeneous, one has
\begin{gather*}
J_{\tau}^{(1/k_1)}\left[\frac{1}{2}\left(1-\frac{x}{\sqrt{b}}\right)\right] = 2^{-p}J_{\tau}^{(1/k_1)}\left(1-\frac{x}{\sqrt{b}}\right),\qquad |\tau| = p.
\end{gather*}
It follows that
\begin{gather*}
\frac{(b)_{\tau}}{(b/2 + (m-1)k_1/2 + (m+3)/4)_{\tau}} J_{\tau}^{(1/k_1)}\left[\frac{1}{2}\left(1-\frac{x}{\sqrt{b}}\right)\right] \approx J_{\tau}^{(1/k_1)}(1)
\end{gather*}
for large positive $b$. Thus, the above Gauss hypergeometric function reduces to
\begin{gather*}
{}_1F_0^{(1/k_1)}\left(\frac{m+1}{2},1\right).
\end{gather*}
Unfortunately, the above series diverges since~\cite{Bak}
\begin{gather*}
{}_1F_0^{(1/k_1)}\left(a,x\right) = \prod_{i=1}^m(1 - x_i)^{-a}.
\end{gather*}
\end{remark}

\subsection[The $B$-type]{The $\boldsymbol{B}$-type}

For this root system, one has
\begin{alignat*}{3}
& R  =  \{\pm e_i, \pm e_i \pm e_j,\  1 \leq i < j \leq m\},\qquad   & &  R_+  =  \{e_i, 1\leq i \leq m, \, e_i \pm e_j,\, 1 \leq i < j \leq m\}, & \\
& S  =  \{e_i - e_{i+1}, \  1 \leq i \leq m, \  e_m \},\qquad & &
C = \{y \in \re^m, \  y_1 > y_2 > \cdots > y_m > 0\}.
\end{alignat*}
The Weyl group is generated by transpositions and sign changes $(x_i \mapsto -x_i)$ and there are two orbits so that $k=(k_0,k_1)$ thereby $\gamma = mk_0 + m(m-1)k_1$ (we assign $k_0$ to $\{\pm e_i, \ 1 \leq i \leq m\}$). The generalized Bessel function\footnote{There is an erroneous sign in one of the arguments in~\cite{Bak}. Moreover, to recover this expression in the $B_m$ case from that given in~\cite{Bak}, one should make the substitutions $a = k_0 - 1/2$, $k_1 = 1/\alpha$, $q = 1 +(m-1)k_1$.} is given by \cite[p.~214]{Bak}
\begin{gather*}
\frac{1}{|W|}D_k^W(x,y) = {}_0F_1^{(1/k_1)}\left(k_0 + (m-1)k_1 + \frac{1}{2}, \frac{x^2}{2t},\frac{y^2}{2t}\right),
\end{gather*}
where  $x^2 = (x_1^2,\dots, x_m^2)$.
Thus
\begin{gather*}
g(x) = |W|\int_C e^{-|y|^2/2}{}_0F_1^{(1/k_1)}\left(k_0 + (m-1)k_1 + \frac{1}{2}, \frac{x^2}{2},\frac{y^2}{2}\right)\prod_{i=1}^m(y_i)V(y^2)dy.
\end{gather*}
The eigenoperator $\mathscr{J}_k$ acts on $W$-invariant functions as
\begin{gather*}
-\mathscr{J}_k^x = \sum_{i=1}^m\partial_i^{2,x} + 2k_0 \sum_{i=1}^m\frac{1}{x_i}\partial_i^x + 2k_1\sum_{i \neq j}\left[\frac{1}{x_i - x_j} + \frac{1}{x_i+x_j}\right] \partial_i^x - E_1^x,
\end{gather*}
therefore $g$ solves the spectral problem with initial value
\begin{gather*}
-\mathscr{J}_k^x g(x) =  m(m+1) g(x),  \qquad g(0) = \int_{\re^m}e^{-|y|^2}\prod_{i=1}^m|y_i| |V(y^2)| dy.
\end{gather*}
 A change of variable $x_i = \sqrt{2y_i}$ shows that $y \mapsto u(y) : = g(\sqrt{2y})$ satisf\/ies
 \begin{gather*}
-\tilde{\mathscr{J}}_k^y u(y)    =   m\frac{(m+1)}{2} u(y),  \qquad u(0) = g(0),\\
-\tilde{\mathscr{J}}_k^y    =   \sum_{i=1}^m y_i\partial_i^{2,y} + 2k_1\sum_{i \neq j} \frac{y_i}{y_i - y_j}\partial_i^y + \left(k_0 + \frac{1}{2}\right)\sum_{i=1}^m\partial_i^y - E_1^y.
\end{gather*}
As a result,
\begin{gather*}
u(y) = u(0){}_1F_1^{(1/k_1)}\left(\frac{m+1}{2}, k_0 + (m-1)k_1 + \frac{1}{2}, y\right),
\end{gather*}
where
\begin{gather*}{}_1F_1^{(1/k_1)}(b,c,z) =  \sum_{p=0}^{\infty}\sum_{\tau}\frac{(b)_{\tau}}{(c)_{\tau}} \frac{J_{\tau}^{(1/k_1)}(z)}{p!}.
\end{gather*}
This can be seen from the dif\/ferential equation (\ref{GF}) and using~\cite{Bak}
\begin{gather*}
\lim_{e \rightarrow \infty}{}_2F_1^{(1/k_1)}\left(e,b,c,\frac{z}{e}\right) = {}_1F_1^{(1/k_1)}(b,c,z).
\end{gather*}
Finally
\begin{gather*}
g\left(\frac{x}{\sqrt t}\right) = g(0){}_1F_1^{(1/k_1)}\left(\frac{m+1}{2}, k_0 + (m-1)k_1 + \frac{1}{2}, \frac{x^2}{2t}\right)
\end{gather*}
and the tail distribution is given by:
\begin{proposition}\label{B1} For $1/2 \leq k_0,k_1 \leq 1$ with either $k_0 > 1/2$ or $k_1 > 1/2$, one has
\begin{gather*}
\p_x^{-l}(T_0 > t) = C_k\prod_{i=1}^m \left(\frac{x_i^2}{2t}\right)^{l_0}\left(V\left(\frac{x^2}{2t}\right)\right)^{2l_1}\\
\phantom{\p_x^{-l}(T_0 > t) =}{}\times e^{-|x|^2/2t}
{}_1F_1^{(1/k_1)}\left(\frac{m+1}{2}, k_0 + (m-1)k_1 + \frac{1}{2} , \frac{x^2}{2t}\right),
\end{gather*}
where $l_0:= k_0 -1/2$, $l_1:= k_1-1/2$ and $V$ stands for the Vandermonde function.
\end{proposition}

\subsection[The $D_m$-type]{The $\boldsymbol{D_m}$-type}

This root system is def\/ined by \cite[p.~42]{Hum}
\begin{gather*}
R  = \{\pm e_i \pm e_j, \ 1 \leq i < j \leq m\},\qquad R_+ = \{e_i \pm e_j, \ 1 \leq i < j \leq m\},
\end{gather*}
and there is one orbit so that $k(\alpha) = k_1$ thereby $\gamma = m(m-1)k_1$. The Weyl chamber is given by
\begin{gather*}
C = \{x \in \re^m,\, x_1 > x_2 > \dots > |x_m|\} = C_1 \cup s_m C_1,
\end{gather*}
where $C_1$ is the Weyl chamber of type $B$ and $s_m$ stands for the ref\/lection with respect to the vector $e_m$ acting by sign change on the variable $x_m$.
\begin{proposition}
For $1/2 < k_1 \leq 1$, the tail distribution writes
\begin{gather*}
\p_x^{-l}(T_0 > t) = C_k \left[V\left(\frac{x^2}{2t}\right)\right]^{2l}e^{-|x|^2/2t}{}_1F_1^{(1/k_1)}\left(\frac{m}{2},(m-1)k_1 + \frac{1}{2}, \frac{x^2}{2}\right).
\end{gather*}
\end{proposition}

\begin{proof} It nearly follows the one given for the root system of type $B$ subject to the following modif\/ications:
if $x \in C_1$, then we perform the change of variable $x_i = \sqrt{2y_i}$, $1 \leq i \leq m$ and for $x \in s_mC_1$ we perform the change of variable $x_i = \sqrt{2y_i}$, $1 \leq i \leq m-1$ and
$x_m = -\sqrt{2y_m}$. In both cases, one gets that $y \mapsto u(y) = g(x^2/2)$ is a symmetric eigenfunction of
\begin{gather}\label{EO}
\sum_{i=1}^m y_i\partial_i^{2,y} + 2k_1\sum_{i \neq j} \frac{y_i}{y_i - y_j}\partial_i^y + \frac{1}{2}\sum_{i=1}^m\partial_i^y - \sum_{i=1}^my_i\partial_i^y
\end{gather}
corresponding to the eigenvalue $m^2/2$. This spectral problem with initial value $1$ at $y=0$ has a unique solution given by
\begin{gather*}
{}_1F_1^{(1/k_1)}(m/2,(m-1)k_1+1/2, y),
\end{gather*}
and the expression of the tail distribution follows.
\end{proof}

\subsection{Second formula}
Suppose that $-1/2 \leq l(\alpha) < 0$ for at least one $\alpha \in R$ and $l(\beta) \geq 0$ for at least one $\beta \in R_+$. Then, by \cite[Proposition~2.15b]{Chy} one writes
\begin{gather*}
\p_x^l ( T_0 > t) = \e_x^0\left[\prod_{\alpha \in R_+}\left(\frac{\langle \alpha,X_t\rangle}{\langle \alpha,x \rangle}\right)^{l(\alpha)}\exp\left(-\frac{1}{2}\sum_{\alpha,\zeta \in R_+}
\int_0^t\frac{\langle \alpha,\gamma\rangle l(\alpha)l(\zeta)}{\langle \alpha,X_s\rangle \langle \zeta,X_s\rangle}ds\right)\right].
\end{gather*}
Using \cite[Proposition~2.15a]{Chy} it follows that
\begin{gather*}
\p_x^l( T_0 > t) = \e_x^r\left[\prod_{\alpha \in R_+}\left(\frac{\langle \alpha,X_t\rangle }{\langle \alpha,x\rangle }\right)^{l(\alpha) - r(\alpha)}
\exp\left(-\frac{1}{2}\sum_{\alpha,\zeta \in R_+}\int_0^t\frac{\langle \alpha,\zeta\rangle l(\alpha,\zeta)}{\langle \alpha,X_s \rangle \langle \zeta,X_s\rangle}ds\right)\right],
\end{gather*}
where $l(\alpha,\zeta) := l(\alpha)l(\zeta) - r(\alpha)r(\zeta)$ and
\begin{gather*}
r(\alpha) = \left\{
\begin{array}{ll}
l(\alpha) & \textrm{if} \  \ l(\alpha) \geq 0,\\
-l(\alpha) & \textrm{if} \ \ l(\alpha) < 0.
\end{array}\right.
\end{gather*}
Then $l(\alpha,\zeta) = 0$ if $l(\alpha)l(\zeta) \geq 0$ and $l(\alpha,\zeta) = -2r(\alpha)r(\zeta)$ otherwise. As a result,
\begin{gather*}
\p_x^l(T_0 > t) = \e_x^r\left[\prod_{\substack{\alpha \in R_+\\l(\alpha) < 0}}\left(\frac{\langle \alpha,x\rangle }{\langle \alpha,X_t\rangle}\right)^{ 2r(\alpha)}
\exp\left(\sum_{\substack{\alpha,\zeta \in R_+\\l(\alpha)l(\zeta) < 0}} \int_0^t\frac{\langle \alpha,\zeta\rangle r(\alpha)r(\zeta)}{\langle \alpha,X_s\rangle \langle \zeta,X_s\rangle}ds\right)\right].
\end{gather*}
Next, note that the exponential functional in the RHS equals $1$ for the irreducible root systems having only one orbit since the sum is empty. Fortunately, this is true for the type $B$.
In fact, note that $l(\alpha)l(\zeta) < 0$ implies that $\alpha$ and $\zeta$ belong to dif\/ferent orbits. Thus, writing
$R_+ = \{e_i, \ 1 \leq i \leq m\} \cup \{e_j \pm e_k,  \ 1 \leq j < k \leq m\}$ so that $\langle e_i,e_j \pm e_k\rangle = \delta_{ij} \pm \delta_{ik}$ gives
\begin{gather*}
 S =  \sum_{i=1}^m\sum_{i < k}\frac{1}{X_t^i}\left[\frac{1}{X_t^i - X_t^k} + \frac{1}{X_t^i + X_t^k}\right]  + \sum_{i=1}^m\sum_{k < i}\frac{1}{X_t^i}\left[\frac{-1}{X_t^k - X_t^i}
+ \frac{1}{X_t^k + X_t^i}\right]
\\
 \phantom{S}{}= \sum_{i=1}^m\sum_{i < k}\frac{2}{(X_t^i)^2 - (X_t^k)^2} - \sum_{i=1}^m\sum_{k < i}\frac{2}{(X_t^k)^2 - (X_t^i)^2} = 0,
\end{gather*}
where $S$ stands for the sum in the above exponential functional (up to a constant). Hence,
\begin{gather*}
\p_x^l(T_0 > t) = \e_x^r\left[\prod_{\substack{\alpha \in R_+\\l(\alpha) < 0}}\left(\frac{\langle \alpha,x\rangle}{\langle \alpha,X_t\rangle }\right)^{2r(\alpha)}\right]
 = \frac{1}{c_kt^{\gamma + m/2}} \prod_{\substack{\alpha \in R_+\\l(\alpha) < 0}} \langle \alpha,x \rangle^{-2l(\alpha)}e^{-|x|^2/2t}
\\
\phantom{\p_x^l(T_0 > t) =}{}\times \int_C e^{-|y|^2/2t}D_k^W\left(\frac{x}{\sqrt t},\frac{y}{\sqrt t}\right)
\prod_{\substack{\alpha \in R_+\\l(\alpha) \geq 0}} \langle \alpha,y\rangle ^{2l(\alpha) +1} \prod_{\substack{\alpha \in R_+\\l(\alpha) < 0}}\langle \alpha,y\rangle dy .
\end{gather*}
For instance, when $k_0 < 1/2$, $k_1 \geq 1/2$ for which $-1/2 \leq l_0 < 0$, $l_1 \geq 0$, the last integral reads:
\begin{align*}
f(x) := \int_C e^{-|y|^2/2t}D_k^W(x,y) V(y^2)^{2l_1 +1} \prod_{i=1}^my_i dy
\end{align*}
while for $k_0 \geq 1/2,\,0 < k_1 < 1/2$, it reads
\begin{align*}
f(x) := \int_C e^{-|y|^2/2t}D_k^W(x,y) V(y^2) \prod_{i=1}^my_i^{2l_0} dy .
\end{align*}
Theorem~\ref{T1} still applies and its corollary needs minor modif\/ications when integrating by parts. The interested reader can easily see that $f$ is an eigenfunction of
$-\mathscr{J}_k$ corresponding to the eigenvalue
\begin{gather*}
m + \sum_{\substack{\alpha \in R_+\\l(\alpha) \geq 0}} [2l(\alpha) + 1] + \sum_{\substack{\alpha \in R_+\\l(\alpha) < 0}} 1 = m + |R_+| +
2\sum_{\substack{\alpha \in R_+\\l(\alpha) \geq 0}} l(\alpha).
\end{gather*}

\section[Tail distribution of $T_0$ and $W$-invariant Dunkl-Hermite polynomials]{Tail distribution of $\boldsymbol{T_0}$ and $\boldsymbol{W}$-invariant\\ Dunkl--Hermite polynomials}

In the sequel, we give an insight to the appearance of the eigenoperator $\mathscr{J}_k$ in our previous computations and show that it is not a mere coincidence. Indeed, this operator already appeared in~\cite{Ros1} and is related to the $W$-invariant counterparts of the so called Dunkl--Hermite polynomials. This fact was behind our attempt to develop an equivalent approach to the pre\-vious one for which the index function is positive. We will give another proof of Theorem~\ref{T1} then express $e^{-|x^2|/2}g(x)$ by means of the $W$-invariant Dunkl--Hermite polynomials. To proceed, we recall some needed facts.

\subsection[Dunkl-Hermite polynomials]{Dunkl--Hermite polynomials}

These polynomials are def\/ined by (see~\cite{Ros1} where they were called generalized Hermite polynomials)
\begin{gather*}
H_{\tau} := e^{-\Delta_k/2}\phi_{\tau}, \qquad \tau = (\tau_1,\dots,\tau_m) \in \mathbb{N}^m,
\end{gather*}
where $(\phi_{\tau})_{\tau \in \mathbb{N}^m}$ are homogeneous polynomials of degree $|\tau| = \tau_1 + \cdots + \tau_m$ and form an orthogonal basis of the vector space of polynomials with real coef\/f\/icients with respect to the pairing inner product introduced in~\cite{Dunkl1}
\begin{gather*}
[p,q]_k = \int_V e^{-\Delta_k/2}p(x)e^{-\Delta_k/2}q(x) \omega_k^2(x) dx
 \end{gather*}
 for two polynomials $p$, $q$ (up to a constant factor). $(H_{\tau})_{\tau}$ are then said to be associated with the basis $(\phi_{\tau})_{\tau}$. By analogy to the one dimensional classical Hermite polynomials, we recall the generating series, the spectral problem and the Mehler-type formula \cite{Dunkl,Ros1}
 \begin{gather}
e^{-|y|^2/2}D_k(x,y)  =  \sum_{\tau}H_{\tau}(x)\phi_{\tau}(y), \label{GS} \\
-\mathscr{J}_k  =  [\Delta_k - \langle x,\nabla\rangle]H_{\tau}(x) =  -|\tau|H_{\tau}(x), \label{SP} \\
\sum_{\tau \in \mathbb{N}^m}H_{\tau}(x)H_{\tau}(y)r^{|\tau|}  =  \frac{1}{(1-r^2)^{\gamma + m/2}}\exp{-\frac{r^2(|x|^2+|y|^2)}{2(1-r^2)}} D_k\left(x,\frac{r}{1-r^2}y\right), \label{MF}
\end{gather}
for $0 < r < 1$.

\subsection[$W$-invariant Dunkl-Hermite polynomials]{$\boldsymbol{W}$-invariant Dunkl--Hermite polynomials}

They are def\/ined up to a constant by
\begin{gather*}
H_{\tau}(x) := \sum_{w \in W} H_{\tau}(wx)
\end{gather*}
and analogs of (\ref{GS}), (\ref{SP}), (\ref{MF}) exist and are derived as follows. Summing twice over $W$ in (\ref{GS}) and using $D_k(wx,w'y) = D_k(x,w^{-1}w'y)$, $w, w' \in W$ \cite{Ros1}
give
\begin{gather*}
\sum_{w,w' \in W} e^{-|y|^2/2}D_k(wx,w'y) = |W| e^{-|y|^2/2}D_k^W(x,y) = \sum_{\tau}H_{\tau}^W(x)\sum_{w \in W}\phi_{\tau}(wy)
\end{gather*}
or equivalently
\begin{gather}\label{I1}
e^{-|y|^2/2}D_k^W(x,y) = \sum_{\tau}H_{\tau}^W(x)\phi_{\tau}^W(y), \qquad \phi_{\tau}^W(y) := \frac{1}{|W|}\sum_{w \in W}\phi_{\tau}(wy).
\end{gather}
Similarly, (\ref{MF}) transforms to
\begin{align}\label{Mehler}
\sum_{\tau \in \mathbb{N}^m}H_{\tau}^W(x)H_{\tau}^W(y)r^{|\tau|} = \frac{|W|}{(1-r^2)^{\gamma + m/2}}\exp{-\frac{r^2(|x|^2+|y|^2)}{2(1-r^2)}} D_k^W\left(x,\frac{r}{1-r^2}y\right).
\end{align}
Finally, summing once over $W$ in (\ref{SP}) and using the fact that $\Delta_k$ and $\langle x,\nabla\rangle$ commute with the action of $W$ (or are $W$-invariant, see \cite[p.~169]{Dunkl}), one gets
\begin{gather}\label{I2}
[\Delta_k - \langle x,\nabla\rangle]H_{\tau}^W(x) = -|\tau|H_{\tau}^W(x).
\end{gather}

 \subsection{Second approach}
 After this wave of formulae, Theorem (\ref{T1}) easily follows after applying (\ref{I2}) to (\ref{I1}) and using the homogeneity of $\phi_{\tau}^W$. Moreover, we derive the following expansion
\begin{proposition}
Let $g$ be the function defined by \eqref{E1}, then
\begin{gather*}
e^{-|x|^2/2}g(x) = \frac{1}{|W|2^{\gamma + m/2}} \sum_{\tau \in \mathbb{N}^m} c_{\tau}^WH_{\tau}^W(x),
\end{gather*}
where
\begin{gather*}
c_{\tau}^W= \frac{1}{2^{|\tau|/2}} \int_Ce^{-|y|^2/2}H_{\tau}^W(y) \prod_{\alpha \in R_+} \langle \alpha,y\rangle dy.
\end{gather*}
\end{proposition}
\begin{proof}
Substituting $D_k^W$ in (\ref{E1}) and using the variable change $y \mapsto ry/(1-r^2)$ yield (since $C$ is a cone)
\begin{align*}
g(x) &= C_{k,r}e^{r^2|x|^2/(2(1-r^2))}\sum_{\tau \in \mathbb{N}^m}H_{\tau}^W(x)r^{|\tau|} \int_Ce^{-|ry|^2/2}H_{\tau}^W\left(\frac{1-r^2}{r}y\right) \prod_{\alpha \in R_+}
\langle \alpha,y\rangle dy
\\ & =C_{k,r}e^{r^2|x|^2/(2(1-r^2))}\sum_{\tau \in \mathbb{N}^m}H_{\tau}^W(x)r^{|\tau|} \int_Ce^{-|y|^2/2}H_{\tau}^W\left(\frac{1-r^2}{r^2}y\right) \prod_{\alpha \in R_+}
\langle \alpha,y\rangle dy \end{align*}
for some constants $C_{k,r}$ depending on $k$, $r$. The result follows by choosing $r = 1/\sqrt 2$.
\end{proof}

Finally, in order to derive $c_{\tau}^W$, one needs

\begin{proposition}[an integration by part formula]
Let $(T_i)_{i=1}^m$ be the Dunkl derivatives~{\rm \cite{Ros1}}. Then, $c_{\tau}^W$ are given by
\begin{gather*}
c_{\tau}^W = (-1)^{|\tau|}\int_C \phi_{\tau}^W(T_1,\dots,T_m)(e^{-|y|^2/2}) \prod_{\alpha \in R_+} \langle \alpha,y\rangle dy.
\end{gather*}
\end{proposition}
\begin{remark}
By Heckman's result (see paragraph at the end of p.~169 in \cite{Dunkl}), $\phi_{\tau}^W(T_1,\dots,T_m)$ acts on $y \mapsto e^{-|y|^2/2}$ as a dif\/ferential operator so that one may perform an integration by parts on the above integral.
\end{remark}
\begin{proof} Recall the Rodriguez-type formula for $H_{\tau}$ \cite{Ros1}
\begin{gather*}
H_{\tau}(y) = (-1)^{|\tau|} e^{|y|^2/2} \phi_{\tau}(T_1,\dots, T_m) \big(e^{-|y|^2/2}\big).
\end{gather*}
The equivariance of the Dunkl operators (see \cite[p.~169]{Dunkl}) yields
\begin{gather*}
H_{\tau}(wy) = (-1)^{|\tau|} e^{|y|^2/2} \phi_{\tau}(\omega^{-1}\cdot)(T_1,\dots, T_m) \big(e^{-|y|^2/2}\big)
\end{gather*}
so that, after summing over $W$, one gets
\begin{gather*}
H_{\tau}^W(y) = (-1)^{|\tau|} e^{|y|^2/2} \phi_{\tau}^W(T_1,\dots, T_m) \big(e^{-|y|^2/2}\big)
\end{gather*}
which ends the proof.
\end{proof}

\section{Examples}

\subsection[$B$-type root systems]{$\boldsymbol{B}$-type root systems} From results in \cite[p.~213]{Bak}, one has on the one hand
\begin{gather*}
H_{2\tau}^W(x) = L_{\tau}^a(x^2/2),\qquad a = k_0-1/2 = l_0
\end{gather*}
and $0$ otherwise, where $L_{\tau}^a$ is the generalized Laguerre polynomial def\/ined in \cite{Lass1} and $\tau$ in the RHS may be seen as a partition of length $m$ since $H_{\tau}^W$ is symmetric (we omitted the dependence on the Jack parameter $1/k_1$ for sake of clarity, see~\cite{Bak} for the details). On the other hand, recall that when
$1/2 \leq k_0, k_1 \leq1$ with either $k_0 > 1/2$ or $k_1 > 1/2$, we derived
\begin{align*}
e^{-|x|^2/2} g(x) & = e^{-|x|^2/2}{}_1F_1^{(1/k_1)}\left(\frac{m+1}{2}, k_0 + (m-1)k_1 + \frac{1}{2} , \frac{x^2}{2}\right)
\\& = {}_1F_1^{(1/k_1)}\left(k_0 - k_1 + m(k_1 - 1/2), k_0 + (m-1)k_1 + \frac{1}{2} , -\frac{x^2}{2}\right)
\end{align*}
by Kummer's relation
\begin{gather*}
e^{-(x_1 + \cdots + x_m)}{}_1F_1^{(1/k_1)}(a,b,x) = {}_1F_1^{(1/k_1)}(b-a,b,-x)
\end{gather*}
whenever it makes sense. Using \cite[Proposition~4.2]{Bak} with $z= (1/2,\dots,1/2)$, $a=k_0 - 1/2$, $q= 1+(m-1)k_1$, $\alpha = 1/k_1$, one gets
\begin{align*}
e^{-|x|^2/2}g(x) &= \frac{1}{2^{(m+1)/2}}\sum_{\tau}\frac{[(l_0 - l_1) + ml_1]_{\tau}}{[k_0 + (m-1)k_1 + 1/2]_{\tau}}
J_{\tau}^{(1/k_1)}\left(\frac{1}{2}\right) L_{\tau}^{k_0-1/2}\left(\frac{x^2}{2}\right)
\\& = \sum_{\tau}\frac{[(l_0 - l_1) + ml_1]_{\tau}}{[k_0 + (m-1)k_1 + 1/2]_{\tau}} \frac{1}{2^{(m+1)/2 + |\tau|}}J_{\tau}^{(1/k_1)}(1)L_{\tau}^{k_0-1/2}\left(\frac{x^2}{2}\right),
\end{align*}
where $l_i = k_i - 1/2$, $i= 1,2$. Here $J_{\tau}^{(1/k_1)}$ is the normalized Jack polynomial denoted $C_{\tau}^{(\alpha)}$ in \cite{Bak} and is related to $\phi_{\tau}^W$ as \cite[p.~201]{Bak}
\begin{gather*}
\phi_{\tau}^W(y) = \frac{(-1)^{|\tau|}}{|\tau|!} \frac{J_{\tau}^{(1/k_1)}(y^2)}{J_{\tau}^{(1/k_1)}(1)}.
\end{gather*}
Comparing the coef\/f\/icients in the expansion of $x \mapsto e^{-|x|^2/2}g(x)$, one gets
\begin{gather*}
\int_C \! J_{\tau}^{1/k_1}(T_1^2,\dots,T_m^2)(e^{-|y|^2/2}) V(y^2)\!\prod_{i=1}^m y_i dy =  \frac{[(k_0 - k_1) + m(k_1 - 1/2)]_{\tau}}{[k_0 + (m-1)k_1 + 1/2]_{\tau}}
\frac{|\tau|![J_{\tau}^{(1/k_1)}(1)]^2}{2^{(m+1+|\tau|)/2}},
\end{gather*}
where $C = \{x_1 > \cdots > x_m > 0\}$.

\subsection[$D$-type root systems]{$\boldsymbol{D}$-type root systems}

Similarly, one  has
\begin{align*}
e^{-|x|^2/2} g(x) &= e^{-|x|^2/2} {}_1F_1^{(1/k_1)}\left(m/2,(m-1)k_1+1/2, \frac{x^2}{2}\right)
\\&=  {}_1F_1^{(1/k_1)}\left((m-1)l_1,(m-1)k_1+1/2, -\frac{x^2}{2}\right)
\\& = \sum_{\tau}\frac{[(m-1)l_1]_{\tau}}{[(m-1)k_1 + 1/2]_{\tau}}L_{\tau}^{-1/2} \left(\frac{x^2}{2}\right) \frac{J_{\tau}^{(1/k_1)}(1)}{2^{m/2 + |\tau|}}.
\end{align*}

\subsection[$A$-type root systems]{$\boldsymbol{A}$-type root systems}

We have already seen that, for $k_1 > 1/2$, $g$ is expressed as a limit of a Gauss multivariate series. With the second approach in hand, $g$ may be expanded by means of the generalized Hermite polynomials def\/ined in \cite{Lass2} and denoted $H_{\tau}$ there, where $\tau$ is a partition of length $m$ (see \cite{Bak} for details).
In particular, this fact may be seen for the value $k=1/2$ as follows, though it has no probabilistic interpretation since $X$ does not hit $\partial C$. Indeed, one seeks an eigenfunction $g$ of
\begin{gather*}
\sum_{i=1}^m \partial_i^{2} + \sum_{i\neq j}\frac{1}{x_i - x_j}\partial_i-  \sum_{i=1}^mx_i \partial_i, \qquad x \in C,
\end{gather*}
associated with the eigenvalue $m+|R_+| = m +m(m-1)/2$ and such that $g(0)= 1$. Easy computations using
\begin{gather*}
\sum_{i \neq j}\frac{x_i}{x_i - x_j} = \frac{m(m-1)}{2}
\end{gather*}
shows that $g(x) = e^{|x|^2/2}$. Using  \cite[Proposition 3.1]{Bak} with $z = (1/\sqrt 2, \dots, 1/\sqrt 2)$ \footnote{We use a rescaling of the generalized Hermite polynomials so that they are orthogonal with respect to $V(x)^{2k} e^{-|x|^2/2}$ rather than $ V(x)^{2k}e^{-|x|^2}$, where $V$ is the Vandermonde function.}, one gets
\begin{align*}
e^{|x|^2/2}  = {}_0F_0^{(2)}(x^2/2)&= e^{1/2}\sum_{\tau} \frac{1}{|\tau|!}H_{\tau}(x^2/2) J_{\tau}^{(1/k_1)}(1/\sqrt 2)
\\& =  e^{1/2}\sum_{\tau} \frac{1}{2^{|\tau|/2}|\tau|!}H_{\tau}(x^2/2) J_{\tau}^{(1/k_1)}(1).
\end{align*}
For general $1/2 < k \leq 1$, one may use the fact that $\phi_{\tau}^W$ is (up to a constant) is a Jack polyno\-mial~$J_{\tau}^{(1/k_1)}$ \cite[p.~190]{Bak} so that
\begin{gather*}
c_{\tau}^W  = (-1)^{\tau} \int_{x_1 > \cdots > x_m} J_{\tau}^{(1/k_1)}(T_1,\dots, T_m)(e^{-|y|^2/2}) V(y) dy
\end{gather*}
where $\tau$ is a partition of length $m$.

\section[Remarks on the first exit times of Brownian motions from $C$]{Remarks on the f\/irst exit times of Brownian motions from $\boldsymbol{C}$}

In \cite{Dou}, authors used a combinatorial approach to write down the tail distribution of the f\/irst exit time from the Weyl chamber of a given root system $R$ by a $m$-dimensional Brownian motion starting at $x \in C$. The result is valid for a wider class of homogeneous $W$-invariant Markov processes and the tail distribution was expressed as Pfaf\/f\/ians of skew symmetric matrices~\cite{Stem}. While we can show, via the heat equation
\begin{gather*}
\frac{1}{2}\Delta u(x,t) = \partial_t u(x,t),\qquad u(x,t) = \mathbb{P}_x^l(T_0 > t),
\end{gather*}
with appropriate boundary values, that our results agree with the ones in \cite{Dou}, we do not succeed to come from Pfaf\/f\/ians to determinants and vice versa. Nevertheless, we do have the following remark: on the one hand, since the $m$-dimensional Brownian motion corresponds to a radial Dunkl process with a null multiplicity function or index function $-l \equiv -1/2$, then our f\/irst approach applies with $k = 1$ so that the Jack parameter $(1/k_1)$ equals one. In this case, it is known that multivariate hypergeometric series of one argument have determinantal representations. For instance, for both $B$ and $D$-types, our formula derived in Subsection~\ref{FF} specialize~to
\begin{gather*}
\p_x^{-1/2}(T_0 > t) =  C\det\left[\left(\frac{x_i^2}{2t}\right)^{m-j+1/2}\fc\left(\frac{m}{2}, m -j + \frac{3}{2}, -\frac{x_i^2}{2t}\right)\right]_{i,j=1}^m, \\
\p_x^{-1/2}(T_0 > t) = C \det\left[\left(\frac{x_i^2}{2t}\right)^{m-j}\fc\left(\frac{m-1}{2}, m -j + \frac{1}{2} , -\frac{x_i^2}{2t}\right)\right]_{i,j=1}^m
\end{gather*}
for both the $B$ and $D$-types root systems  respectively, where $\fc$ is the conf\/luent hypergeometric function. On the other hand, in that cases and for even integers $m$, one has~\cite{Dou}
\begin{gather*}
\p_x^{-1/2}(T_0 > t) = Pf\left(\gamma[(x_i-x_j)/\sqrt{2t}]\gamma[(x_j)/\sqrt{t}]\right)_{1 \leq i,j, \leq m}, \\
\p_x^{-1/2}(T_0 > t) = Pf\left(\gamma[(x_i-x_j)/\sqrt{2t}]\gamma[(x_i + x_j)/\sqrt{2t}]\right)_{1 \leq i,j, \leq m},
\end{gather*}
where
\begin{gather*}
\gamma(a) = \sqrt{\frac{2}{\pi}}\int_0^ae^{-z^2/2} dz,
\end{gather*}
which may be expressed as~\cite{Leb}
\begin{gather*}
 \gamma(a) := \frac{a}{\sqrt{2\pi}}\fc(1/2,  3/2; -a^2/2).
 \end{gather*}
Now, recall that (see \cite[Proposition 2.3]{Stem})
\begin{gather}
\label{PfDet}
{\rm Pf} [\lambda_i\lambda_j]_{1 \leq i < j \leq m} = \prod_{i=1}^m \lambda_i,
 \end{gather}
where in the LHS, the entries of the skew symmetric matrix are $a_{ij} = - a_{ji} = \lambda_i \lambda_j$ for $i < j$. As a result, if $(\lambda_i)_{1 \leq i \leq m}$ are thought of as eigenvalues of some matrix, this formula relates Pfaf\/f\/ians to determinants and we think that it is the way to come from Pfaf\/f\/ians to determinants at least in these cases. We do not have a proof since on the one hand, it is not easy to write down the above determinants as a product of eigenvalues and on the other hand, we are not able to separate the variables $x_i$ and $x_j$ as stated in (\ref{PfDet}).

\subsection*{Acknowledgements}
The author would like to thank C.~Donati Martin for useful remarks and her careful reading of the paper, and P.~Bougerol for explanations of some facts on root systems. He is grateful to M.~Yor for his intensive reading of the manuscript and encouragements.

\pdfbookmark[1]{References}{ref}

\LastPageEnding

\end{document}